\documentclass[11pt]{amsart}
\usepackage{amsmath,amssymb,amsfonts,url}
\numberwithin{equation}{section}

\newtheorem{thm}{Theorem}[section]
\newtheorem{lem}{Lemma}[section]

\newtheorem{prop}{Proposition}[section]

\begin{document}
\title[Construction and non-degeneracy]{Non-degeneracy of Gauss curvature equation with negative conic singularity} 
\author{Juncheng Wei}
\address{Department of Mathematics \\ University of British Columbia \\
 Vancouver, B.C., Canada V6T 1Z2}
\email{jcwei@math.ubc.ca}

\author{Lei Zhang}
\address{Department of Mathematics\\
        University of Florida\\
        358 Little Hall P.O.Box 118105\\
        Gainesville FL 32611-8105}
\email{leizhang@ufl.edu}

\date{\today}

\begin{abstract} We study the Gauss curvature equation with negative singularities. For a local mean field type equation with only one negative index we prove a uniqueness property. For a global equation with one or two negative indexes we prove the non-degeneracy of the linearized equations.
\end{abstract}

\maketitle

\section{Introduction}

In this article we study two closely related equations defined locally and globally in $\mathbb R^2$, respectively. The first equation is define in $\Omega\subset\mathbb R^2$,
which is simply connected, open and bounded.  Throughout the whole article we shall always assume that the boundary of $\Omega$, denoted as $\partial \Omega$, is a rectifiable Jordan curve. We say $\Omega$ is regular from now on.  Let $p_0,p_1,...,p_m\in \Omega$ be a finite set in $\Omega$. Then we consider $v$ as a solution of
\begin{equation}\label{main-eq-1}
\left\{\begin{array}{ll}
\Delta v+\lambda \frac{e^v}{\int_{\Omega}e^v}=-4\pi \alpha_0 \delta_{p_0}+\sum_{i=1}^m 4\pi \alpha_i \delta_{p_i}, \quad \mbox{in }\quad \Omega,  \\
\\
v=0, \quad \mbox{ on }\quad \partial \Omega.
\end{array}
\right.
\end{equation}
where $\alpha_0\in (0,1)$, $\alpha_1,....,\alpha_m>0$ and $\lambda\in \mathbb R$.

The second equation is concerned with the stability of the following global equation: Suppose $u$ is a solution of
\begin{equation}\label{global-e}
\Delta u+e^u=\sum_{i=1}^N4\pi \beta_i \delta_{p_i}, \quad \mbox{ in }\quad \mathbb R^2
\end{equation}
where $\beta_1,...,\beta_n$ are constants greater than $-1$, $p_1,...,p_n$ are the location of singular sources in $\mathbb R^2$. For this equation we shall prove that under some restrictions of $\beta_i$, any bounded solution of the linearized equation has to be the trivial solution.

The background of both equations is incredibly rich not only in mathematics but also in Physics. In particular, the study of (\ref{main-eq-1}) reveals core information on the
configuration of vortices
in the electroweak theory of Glashow-Salam-Weinberg \cite{lai} and Self-Dual Chern-Simons theories \cite{dunne,hong,jackiw}. Also in statistical mechanics the behavior of solutions in
(\ref{main-eq-1}) is closely related to Onsager's model of two-dimensional turbulence with vortex sources \cite{caglioti,chanillo-kiessling}. Most of the motivation and application of both equations come from their connection with conformal geometry. The singular sources represent conic singularities on a surface with constant curvature. There is a large amount of interesting works that discuss the qualitative properties of solutions of such equations. The readers may read into the following works and the references therein \cite{chang,bar-lin-1,bar-lin-2,bar-mal,bar-tan,chanillo-kiessling,chen-lin-wang,chen-lin-10,chen-lin-15,chen-li-duke-93,chen-li-duke-95,li-cmp,lin-wei-ye,luo-tian,
malchiodi-ruiz-1,malchiodi-ruiz-2,nolasco,ohtsuka,spruck-yang,struwe-tarantello,taran-2,taran-1,troyanov-1,troyanov-2,zhang-cmp,zhang-ccm}. It is important to observe that it seems there are very few works to discuss singularities with negative strength and even fewer about the comparison between the negative indexes and positive ones. In this article, using an improved version of the Alexander Bol's inequality we discuss the uniqueness property and the non-degeneracy for a local equation and a global equation. Our proof is based on techniques developed in a number of works of Bartolucci, Lin, Chang-Chen-Lin, etc.

\medskip

To state the main result on the local equation, we first rewrite (\ref{main-eq-1}) using the following Green's function.

For $p\in \Omega$, let $G_{\alpha}(x,p)$ satisfy
$$\left\{\begin{array}{ll}
-\Delta G_{\alpha}(x,p)=4\pi \alpha \delta_p, \quad \mbox{in }\quad \Omega, \\
\\
G_{\alpha}(x,p)=0,\quad x\in \partial \Omega.
\end{array}
\right.
$$
and
$$u=v-G_{\alpha_0}+\sum_{j=1}^mG_{\alpha_j}(x,p_j). $$
Then $u$ satisfies
\begin{equation}\label{main-eq-1-u}
\left\{\begin{array}{ll}
\Delta u+\lambda \frac{He^u}{\int_{\Omega}He^u}=0,\quad \mbox{in }\quad \Omega, \\
\\
u=0,\quad \mbox{on}\quad \partial \Omega.
\end{array}
\right.
\end{equation}
where
\begin{equation}\label{local-h}
H(x)=e^{G_0(x,p_0)}e^{\sum_{j=1}^mG_{\alpha_j}(x,p)}=e^{h(x)}|x-p_0|^{-2\alpha_0}\Pi_{i=1}^m|x-p_i|^{2\alpha_i},
\end{equation}
where $h$ is harmonic in $\Omega$ and is continuous up to the boundary.

The first main result is

\begin{thm}\label{main-thm-1} Let $u$ be a solution of (\ref{main-eq-1-u}) and $H$ be defined by (\ref{local-h}).
 Assume that $\Omega$ is regular, then for any $\lambda\le 8\pi(1-\alpha_0)$ there exists at most one solution to (\ref{main-eq-1}).
\end{thm}

Here we note that for $\lambda<8\pi(1-\alpha_0)$, the existence result has been established by Bartolucci-Malchiodi \cite{bar-mal}. The existence result for $\lambda=8\pi(1-\alpha_0)$ will be discussed in a separate work.

\medskip

The second main result of this article is to consider the nondegeneracy of (\ref{global-e}) when
 there are exactly two negative indexes:
\begin{equation}\label{eq-u-n1}
\left\{\begin{array}{ll}
&\Delta u+e^u=-4\pi \alpha_1\delta_{p_1}-4\pi \alpha_2\delta_{p_2}+\sum_{i=3}^{n} 4\pi\beta_i \delta_{p_i}\quad \mbox{ in }\quad \mathbb R^2, \\
\\
& u(x)=-4\log |x|+\mbox{ a bounded function near } \infty.
\end{array}
\right.
\end{equation}
where $\alpha_1,\alpha_2\in (0,1)$ and $\beta_i>0$ for $i=3,..,n$ and we assume that $n\ge 3$. The assumption of $u$ at infinity says that $\infty$ is not a singularity of $u$ when $\mathbb R^2$ is identified with $\mathbb S^2$.

Let
\begin{equation}\label{eq-u1}
u_1(x)=u(x)+\sum_{i=1}^22\alpha_i\log |x-p_{i}|-2\sum_{i=3}^n \beta_i \log |x-p_i|
\end{equation}
then clearly $u_1$ satisfies
$$\left\{\begin{array}{ll}
\Delta u_1+H_1e^{u_1}=0, \quad \mbox{ in }\quad \mathbb R^2, \\
\\
 u_1(x)=(-4-2\alpha_1-2\alpha_2+2\sum_{i=3}^n \beta_i)\log |x|+O(1),  \mbox{ for } |x|>1,
\end{array}
\right.
$$
where
\begin{equation}\label{v-2}
H_1(x)=\Pi_{i=1}^2|x-p_i|^{-2\alpha_i}\Pi_{i=3}^n|x-p_i|^{2\beta_i}, \quad \mbox{ for }\quad x\in \mathbb R^2.
\end{equation}
Our second main result is
\begin{thm}\label{non-de}Let $u$, $u_1$ and $H_1$ be defined as in (\ref{eq-u-n1}),(\ref{eq-u1}) and (\ref{v-2}), respectively.
Suppose $\phi$ be a classical solution of
\begin{equation}\label{eq-phi}
\Delta \phi+H_1(x)e^{u_1}\phi=0, \quad \mbox{ in }\quad \mathbb R^2.
\end{equation}
If $\lim_{x\to \infty}|\phi(x)|/\log |x|=0$ and $\alpha_1,\alpha_2,\beta_i$ satisfy the following condition:
\begin{equation}\label{index-ineq}
-\max\{\alpha_1,\alpha_2\}+\min\{\alpha_1,\alpha_2\}+\sum_{i=3}^n\beta_i\le 0,
\end{equation}
then $\phi\equiv 0$.
\end{thm}
Here we recall that the total angles at singularities are $2\pi(1-\alpha_1)$, $2\pi(1-\alpha_2)$, $2\pi(1+\beta_i)$ ($i=3,..,n$). For a surface $S$ with conic singularities, let
$$\chi(S,\theta)=\chi(S)+\sum_i(\frac{\theta_i}{2\pi}-1) $$
where $\theta_i$ is the total angle at a conic singularity, $\chi(S)$ is the Euler characteristic of $S$. The purpose of introducing $\chi(S,\theta)$ is to put all surfaces with conic singularities into three cases:
\begin{align*}
&(i) \mbox{ subcritical case if } \quad \chi(S,\theta)<\min_i\{2,\theta_i/\pi\}, \\
&(ii)\mbox{ critical case if } \quad \chi(S,\theta)<\min_i\{2,\theta_i/\pi\}, \\
&(iii)\mbox{ supercritical case if } \quad \chi(S,\theta)>\min_i\{2,\theta_i/\pi\}.
\end{align*}
In our case $\chi(S)=2$ because $S$ is the standard sphere. It is easy to see that (\ref{index-ineq}) refers to the super-critical case. For the subcritical case Troyanov's well known result \cite{troyanov-1} states that every conic singular metric is pointwise conformal to a metric with constant curvature.

Finally if there is only one negative singular source a similar result still holds: Let $u$ satisfy
\begin{equation}\label{eq-u-n2}
\left\{\begin{array}{ll}
\Delta u+e^u=-4\pi \alpha\delta_{p_1}+\sum_{i=2}^{n} 4\pi\beta_i \delta_{p_i}\quad \mbox{ in }\quad \mathbb R^2, \\
\\
 u(x)=-4\log |x|+\mbox{ a bounded function near } \infty.
\end{array}
\right.
\end{equation}
where $\alpha\in (0,1)$ and $\beta_i>0$ for $i=2,..,n$ and we assume that $n\ge 3$.

Let
$$u_1(x)=u(x)+2\alpha\log |x-p_1|-2\sum_{i=2}^n \beta_i \log |x-p_i| $$
then clearly $u_1$ satisfies
\begin{equation}\label{eq-u1-2}
\left\{\begin{array}{ll}
\Delta u_1+H_2e^{u_1}=0, \quad \mbox{ in }\quad \mathbb R^2, \\
 \\
 u_1(x)=(-4-2\alpha+2\sum_{i=2}^n \beta_i)\log |x|+O(1), \quad \mbox{ for }\quad |x|>1,
\end{array}
\right.
\end{equation}
where
\begin{equation}\label{v-e2}
H_2(x)=|x-p_1|^{-2\alpha}\Pi_{i=2}^n|x-p_i|^{2\beta_i}, \quad \mbox{ for }\quad x\in \mathbb R^2.
\end{equation}
Our third main result is
\begin{thm}\label{non-de-2} Let $u_1$ be a solution of (\ref{eq-u1-2}) with $H_2$ defined in (\ref{v-e2}).
Let $\phi$ be a classical solution of
\begin{equation}\label{eq-phi-2}
\Delta \phi+H_2(x)e^{u_1}\phi=0, \quad \mbox{ in }\quad \mathbb R^2.
\end{equation}
If $\lim_{x\to \infty}|\phi(x)|/\log |x|=0$ and $\alpha,\beta_i$ satisfy the following condition:
\begin{equation}\label{index-ineq-2}
-\alpha+\sum_{i=2}^n\beta_i\le 0,
\end{equation}
then $\phi\equiv 0$.
\end{thm}

The organization of this article is as follows. In section two we derive a Bol's inequality with one negative singular source. Then in section three the first two eigenvalues of the linearized local equation is discussed. The proofs of major theorems are arranged in sections 4 and 5. The main approach of this article follows closely from previous works of Bartolucci, Chang, Chen and Lin, etc.

\section{On the Bol's inequality and the first eigenvalues of the local equation}

One of the major tools we shall use is the following Bol's inequality:
\begin{prop} \label{bols-1} Let $\Omega\subset\subset \mathbb R^2$ be a simply connected, open and bounded domain in $\mathbb R^2$. Let $u$ be a solution of
$$\Delta u+Ve^u=0, \quad \mbox{ in }\quad \Omega$$
for
\begin{equation}\label{lem-v}
V=|x-p_1|^{-2\alpha_0}\Pi_{i=2}^n|x-p_i|^{2\beta_i}e^g
\end{equation}
 and $\Delta g\ge 0$ in $\Omega$. Here
$p_1,...,p_n$ ($n\ge 2$) are distinct points in $\Omega$. Let $\omega\subset\Omega$ be an open subset of $\Omega$ such that $\partial \omega$ is a finite union of
 rectifiable Jordan curves. Let
 $$L_{\alpha_0}(\partial \omega)=\int_{\partial \omega}(Ve^u)^{1/2}ds, \quad M_{\alpha_0}(\omega)=\int_{\omega} Ve^u dx. $$
 Then
 \begin{equation}\label{bol-one-neg}
 2L_{\alpha_0}^2(\partial \omega)\ge (8\pi(1-\alpha_0)-M_{\alpha_0}(\omega))M_{\alpha_0}(\omega).
 \end{equation}
 The strict inequality holds if $\omega$ contains more than one singular source or is multiple connected.
\end{prop}

Our proof of Proposition \ref{bols-1} is motivated by the argument in Bartolucci-Castorina \cite{bar-1} and Bartolucci-Lin \cite{bar-lin-1,bar-lin-2}. In fact if $\alpha_0=0$ it was established by Bartolucci-Lin \cite{bar-lin-1}. If $V$ has only singular source at $0$, it was established by Bartolucci-Castorina.  It all starts from an inequality of Huber \cite{huber-annals}:

\emph{ Theorem A (Huber): Let $\omega$ be an open, bounded, simply connected domain with $\partial \omega$ being a rectifiable Jordan curve, $\tilde V=|x|^{-2\alpha_0}e^g$ for some $\Delta g\ge 0$ in $\omega$. Then
$$(\int_{\partial \omega}\tilde V^{1/2}ds)^2\ge 4\pi(1-\alpha_0)\int_{\omega} \tilde Vdx, \quad \mbox{ if }\quad 0\in \omega, $$
$$(\int_{\partial \omega} \tilde V^{1/2} ds)^2\ge 4\pi\int_{\omega} \tilde Vdx, \quad \mbox{ if }\quad 0 \not \in \omega, $$
}

Huber's theorem can be adjusted to the following version

\emph{Theorem B (Bartolucci-Castorina): Let $\omega\subset \mathbb R^2$ be an open bounded domain such that $\partial \omega$ is a rectifiable Jordan curve. Suppose $\bar \omega_B$ is the closure of possibly disconnected bounded component of $\mathbb R^2\setminus \omega$ and $\omega_B$ be the interior of $\bar \omega_B$.Let
$\tilde V=|x|^{-2\alpha_0}e^g$ for some $g$ satisfying $\Delta g\ge 0$ in the interior of $\bar \omega\cup \bar \omega_B$.  Then
$$
(\int_{\partial \omega}\tilde V^{1/2}ds)^2\ge 4\pi(1-\alpha_0)\int_{\omega} \tilde V dx,$$  if $0$ is in the interior of $ \bar \omega\cup \bar \omega_B$.
$$ (\int_{\partial \omega} \tilde V^{1/2} ds)^2\ge 4\pi\int_{\omega} \tilde V dx, $$ if $0$ is not in the interior of $\bar \omega\cup \bar \omega_B$.
}

\medskip

\noindent{\bf Proof of Proposition \ref{bols-1}:} We shall only consider the first case mentioned in Theorem B because the other case corresponds to $\alpha_0=0$.
 Find
$$\left\{\begin{array}{ll}
\Delta q=0,\quad \mbox{in }\quad \omega,\\
q=u,\quad \mbox{ on }\quad \partial \omega.
\end{array}
\right.
$$
and let $\eta=u-q$. Then the equation for $\eta$ is
\begin{equation}\label{2-eq-2}
\left\{\begin{array}{ll}
\Delta \eta+V e^qe^{\eta}=0, \quad \mbox{ in }\quad \omega, \\
\\
\eta=0, \quad \mbox{ on }\quad \partial \omega,
\end{array}
\right.
\end{equation}
and we use
$$
t_m=\max_{\bar \omega} \eta. $$
Then we set
$$\Omega(t)=\{x\in \omega ;  \quad \eta(x)>t\},\quad \Gamma(t)=\partial \Omega(t), \quad \mu(t)=\int_{\Omega(t)}Ve^qdx. $$
Clearly
$\Omega(0)=\omega$, $\mu(0)=\int_{\omega}Ve^qdx$, $\mu(t_m)=\lim_{t\to t_m-}\mu(t)=0$.
Since $\mu$ is continuous and strictly decreasing, it is easy to see that
\begin{equation}\label{mu-d}
\frac{d \mu(t)}{dt}=-\int_{\Gamma(t)}\frac{Ve^q}{|\nabla \eta|} ds, \quad a.e. \quad t\in [0,t_m].
\end{equation}
For all $s\in [0,\mu(0)]$, set
$$\eta^*(s)=|\{t\in [0,t_m],\quad \mu(t)>s\}|$$
where $|E|$ is the Lebesgue measure of the measurable set $E\in \mathbb R$. It is easy to see that $\eta^*$ is the inverse of $\mu$ on $[0,t_m]$ and is continuous, strictly monotone
and differentiable almost everywhere. By (\ref{mu-d}) we have,
for almost all $s\in [0,\mu(0)]$, that
\begin{equation}\label{eta-d}
\frac{d\eta^*}{ds}=-(\int_{\Gamma(\eta^*(s))}\frac{Ve^q}{|\nabla \eta |}dt)^{-1}.
\end{equation}
Let
$$F(s)=\int_{\Omega(\eta^*(s))}e^{\eta}Ve^qdx, \quad a. e. \quad s\in [0,\mu(0)]. $$
Then by the definition of $\Omega(t)$ we see that
$$
F(s)=\int_{\eta^*(s)}^{t_m}e^t(\int_{\Gamma_t}\frac{V e^q}{|\nabla \eta |}ds)dt$$
Using $\beta=\mu(t)$ we further have
\begin{equation}\label{bol-F}
F(s)=\int_0^{s}e^{\eta^*(\beta)}d\beta
\end{equation}
where $\eta^*=\mu^{-1}$ and (\ref{mu-d}) are used.
The definition of $F$ also gives
$$F(0)=\int_{\Omega(\eta^*(0))}e^{\eta }Ve^q=\int_{\Omega(t_m)}e^{\eta}Ve^q=0$$ and
$F(\mu(0))=\int_{\omega}e^{\eta}d\tau=M(\omega)$.
Consequently from (\ref{bol-F}) we obtain
\begin{equation}\label{star-0}
\frac{dF}{ds}=e^{\eta^*(s)},\quad \frac{d^2F}{ds^2}=\frac{d\eta^*}{ds}e^{\eta^*(s)}=\frac{d\eta^*}{ds}\frac{dF}{ds}. \quad a.e.  s.
\end{equation}
Here we use Bartolucci-Castorina's argument in \cite{bar-1} to show that $\eta^*$ is locally lipschitz in $(0,\mu(0))$:
\begin{lem}\label{eta-lip}
For any $0<\bar a\le a<b\le \bar b<\bar u(0)$, there exists $C(\bar a,\bar b,\beta_1,...,\beta_k)>0$ such that
$$\eta^*(a)-\eta^*(b)\le C(b-a).$$
\end{lem}

\noindent{\bf Proof of Lemma \ref{eta-lip}:}
First we find $\Omega_{a,b}$ that satisfies
$$\{x\in \omega;\quad \eta^*(b)\le \eta(x)\le \eta^*(a)\}\subset\subset \Omega_{a,b}\subset\subset \omega. $$
Using Green's representation formula we have
$$|\nabla \eta(x)|\le C+C\int_{\Omega_{a,b}}\frac{1}{|x-y|}|y-p_1|^{-2\alpha_0}dy. $$
Standard estimate gives
\begin{equation}\label{eta-grad}
|\nabla \eta(x)|\le C+C|x-p_0|^{1-2\alpha_0}.
\end{equation}
Recall that $d\eta=Ve^qdx$. Thus
\begin{align*}
&b-a=\mu(\eta^*(b))-\mu(\eta^*(a))\\
=&\int_{\eta>\eta^*(b)}d\tau-\int_{\eta>\eta^*(a)}d\tau\ge \int_{\eta^*(b)<\eta<\eta^*(a)}d\tau\\
=&\int_{\eta^*(b)}^{\eta^*(a)}(\int_{\Gamma(t)}\frac{Ve^q}{|\nabla \eta |}ds)dt
\end{align*}
Using the expression of $V$ in (\ref{lem-v}) and (\ref{eta-grad}) we further have
\begin{align*}
b-a
\ge &\frac 1C\int_{\eta^*(b)}^{\eta^*(a)}\bigg (\int_{\Gamma(t)}\frac{1}{|x-p_0|^{2\alpha_0}+|x-p_0|}\bigg )dt\\
\ge &\frac 1C\int_{\eta^*(b)}^{\eta^*(a)}L_1(\Gamma(t))dt\\
\ge &\min_{\eta^*(b)\le t\le \eta^*(a)}L_1(\Gamma(t))\int_{\eta^*(b)}^{\eta^*(a)}dt\\
\ge &C(\eta^*(a)-\eta^*(b))
\end{align*}
where the estirmate of $\nabla \eta$ was used, $L_1(\Gamma(t))$ stands for the Lebesgue measure of $\Gamma$ and in the last inequality the following standard iso-perimetric inequality $L_1(\Gamma(t))\ge 4\pi |\Omega(t)|\ge 4\pi |\Omega(\eta^*(\bar a)|>0$ is used.
Lemma \ref{eta-lip} is established. $\Box$

\medskip

Now we go back to the proof of Proposition \ref{bols-1}. By Cauchy's inequality

\begin{align}\label{star-1}
(\int_{\Gamma(\eta^*(s))}(Ve^q)^{1/2}ds)^2\le (\int_{\Gamma(\eta^*(s))}\frac{V e^q}{|\nabla \eta|}ds)(\int_{\Gamma(\eta^*(s))}|\nabla \eta |ds) \nonumber\\
=(-\frac{d\eta^*}{ds})^{-1}(\int_{\Gamma(\eta^*(s))}(-\frac{\partial \eta}{\partial \nu})ds), \quad  \mbox{ for } a. e. s\in [0,\mu(0)]
\end{align}
where $\nu=\nabla \eta/|\nabla \eta |$.  Moreover from (\ref{2-eq-2})
\begin{equation}\label{star-2}
\int_{\Gamma(\eta^*(s))}(-\frac{\partial \eta}{\partial \nu})ds=\int_{\Omega(\eta^*(s))}V e^qe^{\eta}dx=F(s), \quad a. e. s\in [0,\mu(0)]
\end{equation}
By Theorem A the following inequality holds for almost all $s\in [0,\mu(0)]$:
\begin{equation}\label{star-3}
(\int_{\Gamma(\eta^*(s))}(Ve^q)^{\frac 12} )^2\ge 4\pi (1-\alpha_0)\mu(\eta^*(s))=4\pi(1-\alpha_0)s.
\end{equation}
Putting (\ref{star-2}) in (\ref{star-1}) yields
\begin{equation}\label{star-4}
(\int_{\Gamma(\eta^*(s))}(Ve^q)^{\frac 12}ds )^2\le (-\frac{d\eta^*}{ds})^{-1}F(s),
\end{equation}

Using (\ref{star-3}) in (\ref{star-4}) we have
$$4\pi(1-\alpha_0)s\le (-\frac{d\eta^*}{ds})^{-1}F(s), \quad a. e. s\in [0,\mu(0)], $$
which is equivalent to
\begin{equation}\label{star-3a}
4\pi(1-\alpha_0)s\frac{d\eta^*}{ds}+F(s)\ge 0, \quad a. e. s\in [0,\mu(0)].
\end{equation}
By (\ref{star-0}) and (\ref{star-3a}), we obtain
$$\frac{d}{ds}[4\pi(1-\alpha_0)(s\frac{dF}{ds}-F(s))+\frac 12 F^2(s)]\ge 0, \quad a. e. s\in [0,\mu(0)]. $$
Let $P(s)$ denote the function in the brackets, then $P$ is well defined, continuous, nondecreasing on $[0,\mu(0)]$. By the Lipschitz property of $\eta^*$, $P$ is absolutely continuous on $[0,\mu(0)]$,
$$P(\mu(0))-P(0)=\lim_{b\to \mu(0)^-}\lim_{a\to 0^+}\int_a^b\frac{dP}{ds}ds.$$
Using $F(0)=0$, $F(\mu(0))=M(\omega)$, $\frac{dF}{ds}|_{s=\mu(0)}=e^0=1$, we have
$$8\pi(1-\alpha_0)(\mu(0)-M(\omega))+M(\omega)^2\ge 0. $$
Then Huber's inequality and $\Gamma(0)=\partial \omega$ further yield
\begin{align*}
2l^2(\partial \omega)&=2(\int_{\partial \omega}(Ve^v)^{1/2}ds)^2\\
&=2(\int_{\partial \omega}(V e^q)^{\frac 12}ds)^2\\
&\ge 8\pi(1-\alpha_0)\mu(0)\\
&\ge M(\omega)(8\pi(1-\alpha_0)-M(\omega)).
\end{align*}
where we have used the fact that $v=q$ on $\partial \omega$.
The Bol's inequality is established. The equality holds if
$Ve^q=|x-p_0|^{-2\alpha_0}|\Phi_t'|^2e^k$ on $\Omega(t)$ for almost all $t\in (0,t_m)$ where $k$ is a constant. In particular for $t=0$ $\Phi_0$ maps $\Omega$ to a ball. In this case $g$ must be harmonic. On the other hand from the equality of Cauchy's inequality we have
$$Ve^q=c_t|\nabla \eta|^2, \quad \mbox{ on }\quad \Gamma(t), \quad a.e.  t\in (0,t_m), $$
for some $c_t>0$. Put $w=\Phi_0(z)$ and $\xi(w)=\eta(\Phi_0^{-1}(w))+k$, we see that $\xi$ satisfies
$$\Delta \xi+|x|^{-2\alpha_0}e^{\xi}=0, $$
and $\xi$ is radial. This $\xi$ is a scaling of
$$\log \frac{8(1-\alpha_0)^2}{1+|x|^{2(1-\alpha_0)})^2}. $$

Thus we have strict inequality in Bol's inequality if at least one of the following situations occurs:
 \begin{enumerate}
 \item $p_1\not \in \omega$,
 \item $\omega$ has at least two singular sources
 \item $\omega$ is not simply connected.
 \end{enumerate}

 $\Box$

\section{The first eigenvalues of the linearized local equation}

\begin{prop}\label{prop-1}
Let $\Omega$ be an open, bounded domain of $\mathbb R^2$ with rectifiable boundary $\partial \Omega$, $V=|x|^{-\alpha_0}\Pi_{i=1}^k|x-p_i|^{2\beta_i}e^g$ for some subharmonic and smooth function $g$, $\alpha_0\in (0,1)$, $\beta_1,...,\beta_k>0$, and we assume that all the singular points are in $\Omega$: $0$, $p_1$,...,$p_k\in \Omega$. Let $w$ be a classical solution of
$$\Delta w+V e^w=0,\quad \mbox{ in }\quad \Omega.$$
Suppose $\hat \nu_1$ is the first eigenvalue of

\begin{equation}\label{star-5}
\left\{\begin{array}{ll}
-\Delta \phi-Ve^w\phi=\hat \nu_1 Ve^w\phi, \quad \mbox{ in }\quad \Omega, \\
\phi=0, \quad \mbox{on}\quad \partial \Omega.
\end{array}
\right.
\end{equation}
Then if $\int_{\Omega}Ve^w\le 4\pi(1-\alpha_0)$  we have $\hat \nu_1> 0$. Moreover if $\int_{\Omega}Ve^w\le 8\pi(1-\alpha_0)$  we have $\hat \nu_2> 0$
\end{prop}

\medskip

\noindent{\bf Proof:} Let $\nu_1=\hat \nu_1+1$ and $\phi$ the eigenfunction corresponding to $\hat \nu_1$, then we have $\phi>0$ and
$$\left\{\begin{array}{ll}
-\Delta \phi=\nu_1Ve^w\phi, \quad \mbox{ in }\quad \Omega, \\
\\
\phi=0, \quad \mbox{ on }\quad \partial \Omega.
\end{array}
\right.
$$
Let
$$U_0(x)=(-2)\log (1+|x|^{2(1-\alpha_0)})+\log (8(1-\alpha_0)^2).$$
Then clearly $U_0$ solves
$$\Delta U_0+|x|^{-2\alpha_0}e^{U_0}=0,\quad \mbox{ in }\quad \mathbb R^2. $$
 For $t\in (0,t_+)$ where $t_+=\max_{\bar \Omega}\phi$, we set
$\Omega(t)=\{x\in \Omega, \quad \phi(x)>t\}$ and we set $R(t)$ to satisfy
$$\int_{\Omega(t)}Ve^w=\int_{B_{R(t)}}e^{U_0}|x|^{-2\alpha_0}. $$
Clearly $\Omega(0)=\Omega$, $R_0=\lim_{t\to 0+}R(t)$, $\lim_{t\to t_+-}R(t)=0$.
Let $\phi^*$ be a radial function from $B_{R_0}\to \mathbb R$. For $y\in B_{R_0}$ and $|y|=r$, set
$$\phi^*(r)=\sup\{t\in (0,t_+)|\quad R(t)>r\}. $$
Then $\phi^*(R_0)=\lim_{r\to R_0-}\phi^*(r)=0$, and the definition implies
$$B_{R(t)}=\{y\in \mathbb R^2, \quad \phi^*(y)>t\}. $$
$$\int_{\phi^*>t}e^{U_0}|x|^{-2\alpha_0}=\int_{\Omega(t)}Ve^w,\quad t\in [0,t_+]. $$
$$\int_{B_{R_0}}|x|^{-2\alpha_0}e^{U_0}|\phi^*|^2=\int_{\Omega}Ve^w\phi^2. $$

Then for almost all $t$
\begin{align}\label{phi-d}
&-\frac{d}{dt}\int_{\Omega(t)}|\nabla \phi|^2=\int_{\phi=t}|\nabla \phi|\\
&\ge (\int_{\phi=t}(Ve^w)^{1/2}ds)^2(\int_{\phi=t}\frac{Ve^w}{|\nabla \phi|}ds)^{-1},\nonumber\\
&=(-\frac{d}{dt}\int_{\Omega(t)}Ve^w)^{-1}(\int_{\phi=t}(Ve^w)^{1/2}ds)^2 \nonumber \\
&\ge \frac 12(8\pi(1-\alpha_0)-\int_{\Omega(t)}Ve^w)(\int_{\Omega_t}Ve^w)(-\frac{d}{dt}\int_{\Omega(t)}Ve^w)^{-1}, \nonumber \\
&=\frac 12(8\pi(1-\alpha_0)-\int_{\phi^*>t}e^{U_0}|x|^{-2\alpha_0})(\int_{\phi^*>t}e^{U_0}|x|^{-2\alpha_0}) \nonumber \\
& \quad \cdot (-\frac{d}{dt}\int_{\phi^*>t}e^{U_0}|x|^{-2\alpha_0})^{-1}. \nonumber
\end{align}
Applying the same computation to $\phi^*$ we see that for almost all $t$, since $\phi^*$ is radial, we have
\begin{align*}
&-\frac{d}{dt}\int_{\Omega(t)}|\nabla \phi^*|^2=\int_{\phi^*=t}|\nabla \phi^*|\\
&=(\int_{\phi^*=t}|x|^{-\alpha_0}e^{U_0/2}ds)^2(\int_{\phi^*=t}\frac{|x|^{-2\alpha_0}e^{U_0}}{|\nabla \phi^*|}ds)^{-1}\\
&=(-\frac{d}{dt}\int_{\Omega(t)}|x|^{-2\alpha_0}e^{U_0})^{-1}(\int_{\phi^*=t}|x|^{-\alpha_0}e^{U_0/2}ds)^2.
\end{align*}
Direct computation on $U_0$ gives
$$(-\frac{d}{dt}\int_{\Omega(t)}|x|^{-2\alpha_0}e^{U_0})^{-1}=\frac 12(8\pi(1-\alpha_0)-\int_{\phi^*>t}e^{U_0}|x|^{-2\alpha_0})(\int_{\phi^*>t}e^{U_0}|x|^{-2\alpha_0}).$$
Thus the combination of the two equations above gives
\begin{align}\label{phi-s-d}
&-\frac{d}{dt}\int_{\Omega(t)}|\nabla \phi^*|^2\\
=&\frac 12(8\pi(1-\alpha_0)-\int_{\phi^*>t}e^{U_0}|x|^{-2\alpha_0})(\int_{\phi^*>t}e^{U_0}|x|^{-2\alpha_0})
(-\frac{d}{dt}\int_{\phi^*>t}e^{U_0}|x|^{-2\alpha_0})^{-1}. \nonumber
\end{align}
for almost all $t\in (0,t_+)$.

Integrating (\ref{phi-d}) and (\ref{phi-s-d}) for $t\in (0,t_+)$ we have
$$\int_{B_{R_0}}|\nabla \phi^*|^2\le \int_{\Omega}|\nabla \phi|^2.$$
If $\nu_1\le 1$, we obtain from (\ref{star-5}) that
\begin{align*}
0\ge (\nu_1-1)\int_{\Omega}Ve^w|\phi |^2=\int_{\Omega} |\nabla \phi |^2-\int_{\Omega}Ve^w|\phi |^2\\
\ge \int_{B_{R_0}}|\nabla \phi^* |^2-\int_{B_{R_0}}e^{U_0}|x|^{-2\alpha_0} |\phi^*|^2.
\end{align*}
Thus the first eigenvalue of
$$-\Delta -|x|^{-2\alpha_0}e^{U_0} $$
on $B_{R_0}$ with Dirichlet boundary condition is non-positive. Since
$$\psi=2(1-\alpha_0)\frac{1-|x|^{2(1-\alpha_0)}}{1+|x|^{2(1-\alpha_0)}}  $$
satisfies
$$-\Delta \psi-|x|^{-2\alpha_0}e^{U_0}\psi=0 \quad \mbox{ in }\quad \mathbb R^2, $$
we see that $R_0\ge 1$. But
$$\int_{B_1}|x|^{-2\alpha_0}e^{U_0}=4\pi(1-\alpha_0),  $$
we clearly have $\hat \nu\ge 0$. From the proof of the Bol's inequality we see that the strictly inequality holds because $\Omega$ has more than one singular points in its interior.

 The proof of $\hat \nu_2>0$ for a higher thresh-hold of $\int_{\Omega}Ve^w$ is very similar. If we consider $\Omega_+$ and $\Omega_-$, which are the set of points where $\phi$ is positive or negative, respectively. Then the integral of $Ve^w$ on at least one of them is less than or equal to $4\pi(1-\alpha_0)$. The argument of re-distribution of mass can be applied to at least one of them. Then we see that either one of them has the integral of $Ve^w$ strictly less than $4\pi(1-\alpha_0)$, which leads to a contradiction, or both regions have their integral equal to
 $4\pi(1-\alpha_0)$. In the latter case the equality cannot hold because $0$ can only be in the interior of at most one region. Then at least one region either does not contain $0$ in its interior, or is not simply connected. The strictly inequality holds in at least one region.  Thus $\hat \nu_2>0$ if $\int_{\Omega}Ve^w\le 8\pi(1-\alpha_0)$.
 $\Box$

\section{Proof of Theorem \ref{non-de}}

First we claim that $\phi$ in the linearized equation is actually bounded.  Recall that $u_1$ satisfies
\begin{align*}
&\Delta u_1+H_1e^{u_1}=0, \quad \mbox{ in }\quad \mathbb R^2, \\
& u_1(x)=(-4+2\alpha_1+2\alpha_2-2\sum_{i=3}^n\beta_i)\log |x|+O(1), \quad \mbox{ at } \quad \infty.
\end{align*}
By the equation for $\phi$ and the mild growth rate of $\phi$ at infinity, we have
$$\phi(x)=\frac{1}{2\pi}\int_{\mathbb R^2}\log |x-y|H_1(y)e^{u_1(y)}\phi(y)dy+c,\quad x\in \mathbb R^2$$
for some $c\in \mathbb R$.

Differentiating the equation above, we have
$$\partial_i\phi(x)=\frac{1}{2\pi}\int_{\mathbb R^2}\frac{x_i-y_i}{|x-y|^2}H_1e^{u_1}\phi(y)dy, \quad i=1,2,\quad x\in \mathbb R^2. $$
By standard estimates in different regions of $\mathbb R^2$, it is easy to see that
$$\partial_i \phi(x)=A\frac{x_i}{|x|^2}+O(|x|^{-1-\delta}),\quad |x|>1 ,\quad i=1,2.$$
for $A=\frac{1}{2\pi}\int_{\mathbb R^2}H_1e^{u_1}\phi$ and some $\delta>0$.
Thus the assumption $\phi(x)=o(\log |x|)$ actually implies
\begin{equation}\label{phi-int-0}
\int_{\mathbb R^2}H_1e^{u_1}\phi=0.
\end{equation}
and
\begin{equation}\label{phi-inf}
\phi(x)=C+O(|x|^{-\delta}),\quad |x|>1
\end{equation}
for some $\delta>0$.

Next we make a transformation on the equation for $u_1$.
Without loss of generality we assume $p_1=0$ and we write $H_1$ as
$$H_1(x)=|x|^{-2\alpha_1}V_1.$$
Let
$$u_2(x)=u_1(\frac{x}{|x|^2})-(4-2\alpha_1)\log |x|, $$
then direct computation shows that
$$\Delta u_2+V_2e^{u_2}=0, \quad \mbox{ in }\quad \mathbb R^2$$
and
$$u_2(x)=(-4+2\alpha_1)\log |x|+O(1),\quad \mbox{ at }\quad \infty. $$
where $V_2(x)=V_1(x/|x|^2)$.
It is also easy to verify that
\begin{equation}\label{int-h1}
\int_{\mathbb R^2}H_1 e^{u_1}=\int_{\mathbb R^2}V_2e^{u_2}.
\end{equation}

Setting $\phi_1(x)=\phi(x/|x|^2)$, we see that
$$\Delta \phi_1+V_2e^{u_2}\phi_1=0, \quad \mbox{ in }\quad \mathbb R^2. $$
Here we note that by the bound of $\phi_1$ near the origin the equation above holds in the whole $\mathbb R^2$.

\medskip

First by the asymptotic behavior of $u_1$ at infinity, integration of the equation for $u_1$ gives
\begin{equation}\label{total-e}
\frac{1}{2\pi}\int_{\mathbb R^2}H_1e^{u_1}=4-2(\alpha_1+\alpha_2)+2\sum_{i=3}^n\beta_i\le 4(1-\alpha_2).
\end{equation}

From the definition of $\phi$ we have
$\phi_1(x)\to c_0$ as $x\to \infty$ for some $c_0\in \mathbb R$.
Without loss of generality we assume $c_0\le 0$. By the same estimate for $\phi$ we have
\begin{equation}\label{inte-0}
\int_{\mathbb R^2}V_2e^{u_2}\phi_1=0.
\end{equation}

By (\ref{int-h1}) and (\ref{total-e}) we have
$$\int_{\mathbb R^2}V_2e^{u_2}\le 8\pi (1-\alpha_2). $$

Let $\phi_2$ be an eigenfunction corresponding to eigenvalue $\hat \nu$:
$$\left\{\begin{array}{ll}
-\Delta \phi_2-V_2e^{u_2}\phi_2=\hat \nu V_2e^{u_2}\phi_2, \quad \mbox{in}\quad \mathbb R^2, \\
\lim_{x\to\infty}\phi_2(x)=c_0\le 0, \\
\int_{\mathbb R^2}V_2e^{u_2}\phi_2=0.
\end{array}
\right.
$$
We claim that $\hat \nu>0$.

By way of contradiction we assume that $\hat \nu\le 0$. By setting $\nu=1+\hat \nu$ we clearly have $\nu\le 1$ and
$$\Delta \phi_2+\nu V_2e^{u_2}\phi_2=0,\quad \mbox{in}\quad \mathbb R^2. $$

Let $\Omega^+=\{x;\quad \phi_2(x)>c_0\}$ then by the same argument as in the proof of the previous proposition we must have
$$\int_{\Omega^+}V_2e^{u_2}=c_2(c_0)\ge 4\pi(1-\alpha_2) $$
and if the equality holds, we have $c_0=0$, there is one singular source with negative index $-4\pi \alpha_2$ in the interior of $\Omega_+$, which has to be simply connected at the same time. All other singular sources (which have positive indexes) are not in the interior of $\Omega_+$.

 Let $\phi^*$ be the rearrangement of $\phi_2$ in $\Omega_+$. By the previous argument we have
$$\int_{\Omega_+}|\nabla \phi_2|^2\le \int_{B_{R_1}}|\nabla \phi^*|^2 $$ and
$c_2(c_0)=\int_{B_{R_1}}|x|^{-2\alpha_2}e^{U_0}$.
Let
$$c_1=\min_{\mathbb R^2}\phi_2 $$
and we set $R_2$ to make
$$\int_{B_{R_2}\setminus B_{R_1}}|x|^{-2\alpha_2}e^{U_0}=\int_{\mathbb R^2\setminus \Omega_+}V_1e^{u_2}. $$
Note that $R_2$ could be $\infty$.
Then we define a radial function $\phi^{**}$ from $B_{R_2}\setminus B_{R_1}\to \mathbb R$:
for any $y\in B_{R_2}\setminus B_{R_1}$, $|y|=r$,
$$\phi^{**}(r)=\inf\{t\in (c_1,c_0)|\quad R^{(-)}(t)<r\}, $$
where $R^{(-)}(t)$ is defined by
$$\int_{B_{R_2}\setminus B_{R^{(-)}(t)}}|x|^{-2\alpha_2}e^{U_0}=\int_{\phi_2<t}V_2e^{u_2},\quad \forall t\in (c_1,c_0). $$
The definition of $\phi^{**}$ implies
$$\int_{B_{R_2}\setminus B_{R_1}}|x|^{-2\alpha_2}e^{U_0}|\phi^{(**)}|^2=\int_{\Omega^-}V_2e^{u_2}|\phi_2 |^2, \quad \Omega^-=\mathbb R^2\setminus \Omega_+, $$
and
$$\int_{B_{R_2}\setminus B_{R_1}}|x|^{-2\alpha_2}e^{U_0}\phi^{(**)}=\int_{\Omega^-}V_2e^{u_2}\phi_2 , \quad \Omega^-=\mathbb R^2\setminus \Omega_+. $$
The symmetrization also gives
$$\int_{B_{R_2}\setminus B_{R_1}}|\nabla \phi^{**}|^2\le \int_{\Omega^-}|\nabla \phi_2 |^2. $$
Now we set
$$\phi_*:B_{R_2}\to \mathbb R, \quad \phi_* \mbox{ radial }
\phi_*(r)=\left\{\begin{array}{ll}
\phi^*(r), \quad r\in [0, R_1], \\
\\
\phi^{**}(r),\quad r\in [R_1, R_2).
\end{array}
\right.
$$
Since $\phi_*$ is continuous, monotone, we have
$$
\int_{B_{R_2}}|\nabla \phi_*|^2\le \int_{\mathbb R^2}|\nabla \phi_2|^2=\int_{\mathbb R^2}V_2e^{u_2}|\phi_2 |^2=\int_{B_{R_2}}|x|^{-2\alpha_2}e^{U_0}|\phi_*|^2.
$$
From the definition of $\phi_*$ we also have
$$\int_{B_{R_2}}|x|^{-2\alpha_2}e^{U_0}\phi_*=0. $$

Let
\begin{align*}
K^*=&\inf \{\int_{\mathbb R^2}|\nabla \psi|^2dx, \quad \psi \mbox{ is radial}, \\
& \int_{\mathbb R^2}|x|^{-2\alpha_2}e^{U_0}\psi=0,\quad
\int_{\mathbb R^2}|x|^{-2\alpha_2}e^{U_0}\psi^2=1.\}.
\end{align*}

By H\"older's inequality we have
$$
|\int_{\mathbb R^2}|x|^{-2\alpha_2}e^{U_0}\psi dx|\le (\int_{\mathbb R^2}|x|^{-2\alpha_2}e^{U_0}\psi^2)^{\frac 12}(\int_{\mathbb R^2}|x|^{-2\alpha_2}e^{U_0})^{1/2}.
$$
Which implies that the minimizer (say $\psi^*$) also satisfies
$$\int_{\mathbb R^2}|x|^{-2\alpha_2}e^{U_0}\psi^*=0. $$
Clearly the minimizer $\psi^*$ satisfies
$$\Delta \psi^*+K^*|x|^{-2\alpha_2}e^{U_0}\psi^*=0,\quad \mbox{ in }\quad \mathbb R^2. $$
From $\phi_*$ and the definition of $K^*$ we already know $K^*\in (0,1)$.
Our goal is to show that $K^*=1$ by an argument of Chang-Chen-Lin \cite{chang}.  $\psi^*$ should only change sign once. Let $\xi_0$ be the zero of $\psi^*$.

Integrating the equation for $\psi^*$, we have
$$r\frac{d}{dr}\psi^*(r)=-K^*\int_0^r|s|^{1-2\alpha_2}e^{U_0(s)}\psi^*(s)ds=K^*\int_r^{\infty}s^{1-2\alpha_2}e^{U_0}\psi^*(s)ds<0 $$
for $r>\xi_0$. Thus $\psi^*$ is decreasing for $r\ge \xi_0$ and $r\frac{d}{dr}\psi^*(r)\to 0$ as $r\to \infty$. The equation for $\psi^*$ also gives
$$|r\frac{d}{dr}\psi^*(r)|\le K^*(\int_r^{\infty}|s|^{1-2\alpha_2}e^{U_0}(\psi^*(s))^2ds)^{1/2}(\int_r^{\infty}s^{1-2\alpha_2}e^{U_0(s)}ds)^{\frac 12}\le Cr^{-1} $$
for large $r$. Therefore $\lim_{r\to \infty}\psi^*(r)$ exists and is a negative constant.

Let
$$\psi(r)=2(1-\alpha_2)\frac{1-r^{2(1-\alpha_2)}}{1+r^{2(1-\alpha_2)}}.$$
Then $\psi$ satisfies
$$\Delta \psi+r^{-2\alpha_2}e^{U_0}\psi=0, \quad \mbox{ in }\quad \mathbb R^2.
$$
It is easy to obtain from the equation for $\psi$ and the one for $\psi^*$ the following:
$$r(\frac{\psi^*}{\psi(r)})'=\frac{1-K^*}{\psi^2(r)}\int_0^rs^{1-2\alpha_2}e^{U_0(s)}\psi^*(s)\psi(s)ds. $$
 If $\xi_0<1$, $\frac{\psi^*(r)}{\psi(r)}$ is increasing from $r\in (0,\xi_0]$. Clearly this is not possible because otherwise this could happen:
$$0<\frac{\psi^*(0)}{\psi(0)}<\frac{\psi^*(\xi_0)}{\psi(\xi_0)}=0. $$
On other hand we observe that it is also absurd to have $\xi_0>1$, indeed, had this happened, we would start from
$$\lim_{R\to \infty}R(\frac{\psi^*}{\psi})'(R)\psi^2(R)-r(\frac{\psi^*}{\psi})'(r)\psi^2(r)=(1-K^*)\int_r^{\infty}s^{1-2\alpha_2}e^{U_0}\psi^*(s)\psi(s)ds. $$
Since
$$\lim_{R\to \infty} R(\frac{d}{dr}\psi^*(R)\psi(R)-\psi'(R)\psi^*(R))=0, $$
we have
$$-r(\frac{\psi^*}{\psi})'\psi^2(r)=(1-K^*)\int_r^{\infty}s^{1-2\alpha_2}e^{U_0(s)}\psi^*(s)\psi(s)ds. $$
If $\xi_0>1$, $\frac{\psi^*(r)}{\psi(r)}$ is decreasing for $r>1$, which yields
$$0=\frac{\psi^*(\xi_0)}{\psi(\xi_0)}>\lim_{r\to \infty}\frac{\psi^*(r)}{\psi(r)}=-\frac{1}{2(1-\alpha_2)}\lim_{r\to \infty}\psi^*(r)>0. $$
This contradiction proves that $\xi_0=1$ and $\psi^*(r)\psi(r)>0$ for all $r\neq 1$. Furthermore
\begin{align*}
0&=\lim_{r\to \infty} (\frac{d}{dr}\psi^*(r)\psi(r)-\frac{d}{dr}\psi(r)\psi^*(r))r\\
&=(1-K^*)\int_0^{\infty}s^{1-2\alpha_2}e^{U_0}\psi^*(s)\psi(s)ds.
\end{align*}
Thus we have proved that $K^*=1$ and the desired contradiction. Theorem \ref{non-de} is established. $\Box$

\medskip

The proof of Theorem \ref{non-de-2} is very similar, we just use Kelvin transformation to move the negative singularity to infinity, then use the same argument with the standard Bol's inequality for non-negative indexes.

\section{The proof of Theorem \ref{main-thm-1}.}
 Our argument follows from a previous result of Bartolucci-Lin \cite{bar-lin-2} for non-negative indexed singularities. We prove by way of contradiction. Suppose $u$ is a solution of (\ref{main-eq-1-u})
 and a nonzero function $\tilde \phi\in H^1_0(\Omega)$ is a solution of
$$\left\{\begin{array}{ll}
-\Delta \tilde \phi-\lambda \frac{He^u}{\int_{\Omega}He^udx}\tilde \phi+\lambda(\int_{\Omega}He^u\tilde \phi)\frac{He^u}{(\int_{\Omega}He^u)^2}=0, \mbox{ in } \Omega, \\
\\
\tilde \phi=0, \quad \mbox{ on }\quad \partial \Omega.
\end{array}
\right.
$$
Let $w=u+\log \lambda-\log (\int_{\Omega}He^udx)$ and
$$\phi=\tilde \phi-\frac{\int_{\Omega}He^u\tilde \phi}{\int_{\Omega}He^u}, $$
we have
\begin{equation}\label{phi-th1}
\left\{\begin{array}{ll}
\Delta \phi+He^w\phi=0, \quad \mbox{ in }\quad \Omega, \\
\phi=c_0, \quad \mbox{ on }\quad \partial \Omega, \\
\int_{\Omega}He^w\phi=0, \\
\lambda=\int_{\Omega}He^w\le 8\pi(1-\alpha_0).
\end{array}
\right.
\end{equation}
Without loss of generality we assume $c_0\le 0$. Our goal is to show that $\phi\equiv c_0$, which further leads to $c_0=0$ obviously. If $c_0=0$, $\phi$ must change sign if not identically equal to $0$. But this situation is ruled out by Proposition \ref{prop-1} that $\nu_2>0$. So we only consider $c_0<0$. Let
$$\Omega_+=\{x\in \Omega, \quad \phi(x)>0\quad \}, \quad \Omega_-=\{x\in \Omega, \quad \phi(x)<0\quad \}. $$
Clearly $dist(\Omega_+,\partial \Omega)>0$ Then if $\int_{\Omega_+}He^w\le 4\pi(1-\alpha_0)$ there is no way for $\phi$ to satisfy (\ref{phi-th1}) on $\Omega_+$ without being identically zero. Then using the same rearrangement argument as in the proof of Theorem \ref{non-de} we can also reach the following conclusion:
If $\phi_2$ is a solution of
$$\left\{\begin{array}{ll}
-\Delta \phi_2-\lambda e^uw\phi_2=\nu e^u w \phi_2, \quad \mbox{ in } \quad \Omega, \\
\phi_2=c_0, \quad \mbox{ on }\quad \partial \Omega,
\end{array}
\right.
$$
then $\nu>0$. The remaining part of the the proof of Theorem \ref{main-thm-1} follows by standard argument in \cite{chang} and \cite{bar-lin-2}. We include it with necessary modification.

If we use $L_{\lambda}$ to denote the linearized operator of (\ref{main-eq-1-u}), we know that all eigenvalues of  $L_{\lambda}$ are strictly positive for $\lambda\in [0, 8\pi(1-\alpha_0)]$. By using the improved Moser-Trudinger inequality \cite{malchiodi-ruiz-1} one can easily find a solution of (\ref{main-eq-1-u}) by direct minimization method. By the uniform estimate of the linearized equation and standard elliptic estimate we have: for any $\epsilon \in (0, 8\pi(1-\alpha_0))$,
\begin{equation}\label{unif-es}
\|u_{\lambda}\|_{\infty}\le \lambda C_{\epsilon}
\end{equation}
for some $C_{\epsilon}>0$,  $\lambda\in [0, 8\pi(1-\alpha_0)]$ and $u_{\lambda}$ as solution of (\ref{main-eq-1-u}). Let $S_{\lambda}$ be the solution's branch for (\ref{main-eq-1-u}) bifurcating from
$(u,\lambda)=(0,0)$. The standard bifurcation theory of Crandall-Rabinowitz \cite{crandall} gives that $S_{\lambda}$ is a simple branch near $\lambda=0$. Which means for $\lambda>0$ small there exists one and only solution for (\ref{main-eq-1-u}) and $S_{\lambda}$ is smooth in $C^2(\Omega)\times \mathbb R$. By the implicit function theorem (because $L_{\lambda}$ has positive first eigenvalue) $S_{\lambda}$ can be extended uniquely for $\lambda\in (0, 8\pi(1-\alpha_0))$. If for any given $\lambda\in (0, 8\pi(1-\alpha_0)$ there is another
solution, it implies the other solution's branch
 does not bend in $[0, 8\pi(1-\alpha_0))$. By the uniform estimate (\ref{unif-es}) this second branch intersects $S_{\lambda}$ at $(u,\lambda)=(0,0)$. This contradiction proves the uniqueness for $\lambda\in [0, 8\pi(1-\alpha_0))$. If a solution exists for $\lambda=8\pi(1-\alpha_0)$, the implicit function theorem and the uniqueness result can be combined to prove the uniqueness in this case as well.
Theorem \ref{main-thm-1} is established. $\Box$

\end{document}